\documentclass{article}

\usepackage{amsmath,amssymb,makeidx,amsfonts,latexsym,theorem,enumerate}

\newtheorem{lemma}{Lemma}[section]
\newtheorem{theorem}{Theorem}[section]
\newtheorem{corollary}{Corollary}[section]

\begin{document}

\author{Matthias Gundlach
\\ Institut f\"ur Dynamische Systeme, Universit\"at Bremen\\
Postfach 330 440, 28334 Bremen, Germany
\and Andrei Khrennikov\footnote{This research was supported by the
  visiting professor fellowship at Bremen University.} \& Karl-Olof Lindahl 
\footnote{This research was supported by the grant
''Strategical investigations'' of V\"axj\"o University.}\\
School of Mathematics and Systems Engineering\\
V\"{a}xj\"{o} University, 351 95, V\"{a}xj\"{o}, Sweden}
\title{On ergodic behavior of $p$-adic dynamical systems
\footnote{This note was published in \emph{Infinite Dimensional Analysis}, Vol. 4,
No. 4  (2001) 569--577.}}

\date{\empty}
\maketitle

\begin{abstract}
Monomial mappings, $x\mapsto x^n$, are topologically transitive and
ergodic with respect to Haar measure on the unit circle in the
complex plane.
In this paper we obtain an anologous result for monomial dynamical systems
over $p-$adic
numbers. The process is,
however, not straightforward. The result will depend on the natural number
$n$. Moreover,
in the $p-$adic case we never have ergodicity on the unit circle, but on
the circles
around the point $1$.  
\end{abstract}

\section{Introduction}

Investgations in $p-$adic quantum physics 
  \cite{VV} -- \cite{Kh1} (especially string theory   \cite{VV},  
  \cite{Freund}, \cite{Vladimirov1})
stimulated an increasing interest in studying $p-$adic dynamical systems,
see for example \cite{Thiran et al} -- \cite{Bosio}.
Some steps
in this direction \cite{Kh1}, \cite{Albeverio3} 
demostrated that even the simplest (monomial) discrete dynamical
systems over the fields of $p-$adic numbers $\mathbb{Q}_p$ have quite complex
behavior. This behavior depends crucially on the prime number $p>1$ (which
determines
$\mathbb{Q}_p$). By varying $p$ we can transform attractors into centers of
Siegel discs
and vice versa. The number of cycles and their lengths also depend
crucially on $p$ \cite{Pezeda}, \cite{Nilsson}. Other aspects (such as recurrence,
renormalization and conjugacy)
of monomials and the dynamics of the endomorphisms of the multiplicative group in general was
investigated by Arrowsmith and Vivaldi \cite{Arrowsmith Vivaldi2}. 

Some applications of discrete $p-$adic dynamical systems to cognitive
sciences
and neural networks were considered in \cite{Kh1}, 
\cite{Albeverio4}. Some of these cognitive models are described
by random dynamical systems in the fields $p-$adic numbers, see
\cite{Dubischar}, \cite{Li1}. 
In this note we study ergodicity of monomial $p-$adic dynamical
systems on spheres.
For a system $\psi_n(x)=x^n$, $n=2,3,...,$ the result depends crucially on
the relation
between $n$ and $p$. Our proof is essentialy based on $p-$adic analysis
(analytic mappings),
\cite{Schikhof}. We remark that the corresponding fact for the field of
complex numbers,
$\psi_n: \mathbb{C}\rightarrow \mathbb{C}$, for the sphere
$\left|z\right|=1$, is rather
trivial, see \cite{Walters}.

\section{$p$-adic numbers}\label{p-adic numbers}

The system of $p$-adic numbers $\mathbb{Q}_p$ was constructed by K. Hensel 
in the 1890s.

The field of real numbers $\mathbb{R}$ is constructed as the completion of the
field of rational
numbers $\mathbb{Q}_p$ with respect to the metric $\rho(x,y)$ $=$ $\vert x - y
\vert$, where $\vert
\cdot\vert$ is the usual valuation of distance given by the absolute value.
The fields of $p$-adic
numbers $\mathbb{Q}_p$  are constructed in a corresponding way, but  using
other  valuations.
For a prime number $p$, the $p$-adic valuation $\vert \cdot\vert_p $ is
defined in the following
way. First we define it for natural numbers. Every natural number $n$ can
be represented as the
product of prime  numbers, $n$ $=$ $2^{r_2}3^{r_3} \cdots p^{r_p} \cdots$,
and we define
$\vert n\vert_p$ $=$ $p^{-r_p}$, writing $\vert 0 \vert_p$ $=0$  and
$\vert -n\vert_p$ $=$
$\vert n \vert_p$.  We then extend the definition of the $p$-adic valuation
$\vert\cdot\vert_p$
to all rational numbers by setting 
\[
\vert n/m\vert_p= \vert n\vert_p/\vert m\vert_p
\]
for $m$ $\not=$ $0$. The completion of $\mathbb{Q}_p$ with respect to the 
metric
$\rho_p (x,y)$ $=$
$\vert x- y\vert_p$ is the locally compact field of $p$-adic numbers 
$\mathbb{Q}_p$.

The number fields $\mathbb{R}$ and $\mathbb{Q}_p$ are unique in a
 sense, since by Ostrovsky's
theorem (\cite{Schikhof}), $\vert \cdot \vert$ and  $\vert\cdot\vert_p$ are the
only possible  valuations
on $\mathbb{Q}$, but have quite distinctive properties. The field of real
numbers $\mathbb{R}$ with its
usual valuation satisfies $\vert n\vert$ $=$ $n$ $\to$ $\infty$ for
valuations of natural numbers $n$
and is said to be {\it Archimedean.\/} By a well know theorem of number
theory \cite{Schikhof}, the only
complete Archimedean fields are those of the real and the complex numbers.
In contrast, the fields
of $p$-adic numbers, which satisfy $\vert n\vert_p$ $\leq$ $1$ for all $n$
$\in\mathbb{N}$, are
examples of {\it non-Archimedean\/} fields.

Unlike the absolute value distance $\vert \cdot \vert$, the $p$-adic
valuation satisfies the strong
triangle inequality
\begin{equation}\label{str}
|x+y|_p \leq  \max[|x|_p,|y|_p],  \quad x,y \in \mathbb{Q}_p,
\end{equation}
with equality in the case that $|x|_p\neq |y|_p$.

Write $B_r(a)$ $=$ $\{x\in \mathbb{Q}_p: |x -a|_p \leq r\}$ and $B_r^-(a)$ $=$
$\{x\in \mathbb{Q}_p:
|x -a|_p < r\}$ where $r$ $=$ $p^n$ and $n$ $=$ $0$, $\pm 1$, $\pm 2$,
$\ldots$.
These are the ``closed'' and ``open''  balls in $\mathbb{Q}_p$ while the sets
$S _r(a)$ $=$ $\{x \in K: |x -a|_p = r \}$ are the spheres in $\mathbb{Q}_p$
of such radii $r$.
Any $p$-adic
ball $B_r(0)$ is an additive subgroup of $\mathbb{Q}_p$, while the ball
$B_1(0)$ is also
a ring, which is called the {\it ring of $p$-adic integers} and is denoted
by $\mathbb{Z}_p$.

As a consequence of the strong triangle property (\ref{str}) every element
$b$ of a ball $B_r(a)$ is the center of this
ball in the sense that $B_r(b)=B_r(a)$ for every $b\in B_r(a)$.

Each $p-$adic number has a unique expansion of the form
$\sum_{n=M}^{\infty}a_np^n$, where $M$ is an integer and  $a_n\in \{0, ...,p-1\}$.
The $p-$adic integers $\mathbb{Z}_p$, are $p-$adic
numbers with expansions $\sum_{n=0}^{\infty}a_np^n$. 
The absolute value of a $p-$adic number is $p^{-k}$ if and only if $k$ is the 
first
digit in the expansion which differs from zero. Thus 
$\left|x-y\right|_p\leq p^{-k}$ if and only if
$x\equiv y\mod p^k$.

Let $(m,n)$ denote the greatest common divisor of $m$ and $n$.  
We say that $n$ is a (multplicative) \emph{unit} (with respect to the
prime number $p$) iff $(n,p)=1$. Let $G_{p^l}$, $l\geq 1$, be the
multiplicative group of units  in the residue field modulo $p^l$.
Let us by $\langle n\rangle =\{n^N:\quad N\in\mathbb{N}\}$ denote the set
\emph{generated} by $n$.

Let $\psi_n$ be a (monomial) mapping on $\mathbb{Z}_p$ taking $x$ to $x^n$. 
Then all spheres
$S_{p^{-l}}(1)$ are 
$\psi _n$-invariant iff $n$ is a multplicative unit. 
This is a consequence of the following result in $p-$adic analysis.
A proof is given in Appendix \ref{lemma}. 

\begin{lemma}
\label{lemma 1}Let $x,y\in S_1(0)$ and suppose $\left| x-y\right| _p<1.$
Then for all natural numbers $n,$%
\begin{equation}
\left| x^n-y^n\right| _p\leq \left| n\right| _p\left| x-y\right| _p,
\label{fundamental-inequality}
\end{equation}
with equality for $p>2$. Moreover equality also holds for $p=2$ if $n$
is odd.
\end{lemma}

\vspace{1.5ex}
\noindent
In particular, $\psi_n$ is an isometry on $S_{p^{-l}}(1)$ if and only if $(n,p)=1$.
Therefore we will henceforth
assume that $n$ is a
unit.  Also note that, as a consequence, $S_{p^{-l}}(1)$ is not a group under 
multiplication. Thus our investigations are not about the dynamics on a 
compact (abelian) group.

\section{Minimality}

Let us consider the dynamical system $x\mapsto x^n$ on spheres 
$S_{p^{-l}}(1)$. 
The result depends crucially on the following well known result from group 
theory.
 
\begin{lemma}\label{generation}
Let $p>2$ and $l$ be any natural number, then the natural number $n$ is a 
generator of
$G_{p^l}$ if and only if $n$ is a generator of $G_{p^2}$. $G_{2^l}$ is 
noncyclic for
$l\geq 3$.
\end{lemma}

\vspace{1.5ex}
\noindent
Recall that a dynamical system given by a continuous transformation $\psi$ on
a compact metric space $X$ is called \emph{topologically transitive} if there 
exists
a dense orbit $\{\psi^n(x): n\in\mathbb{N}\}$ in $X$, and (one-sided)
\emph{minimal}, if all orbits for $\psi$ in $X$ are dense. For the case of 
monomial systems $x\mapsto x^n$ on spheres $S_{p^{-l}}(1)$ topological 
transitivity
means the existens of an $x\in S_{p^{-l}}(1)$ s.t. each $y\in S_{p^{-l}}(1)$
is a limit point in the orbit of $x$, i.e. can be represented as
\begin{equation}\label{limit point}
y=\lim_{k\rightarrow\infty} x^{n^{N_k}},
\end{equation}
for some sequence $\{N_k\}$, while minimality means that such a property holds
for any $x\in S_{p^{-l}}(1)$. Our investigations are based on the following 
theorem.

  \begin{theorem}
    For $p\neq 2$ the set $\langle n\rangle$ is dense in
    $S_1(0)$ if and only if $n$ is a generator of $G_{p^2}$. 
  \end{theorem}
    {\noindent \bf Proof. }%
      We have to show that for every $\epsilon >0$ and 
      every $x\in S_1(0)$ there is a $y\in \langle n\rangle$ such 
      that $\left|x-y\right|_p<\epsilon$ if and only if $n$ is a generator of $G_{p^2}$. 
      First, assume that $n$ is a generator of $G_{p^2}$.
      Let $\epsilon>0$ and $x\in S_1(0)$ be
      arbitrary. Because of the discreteness of the $p-$adic metric we can assume that
      $\epsilon=p^{-k}$ for some natural number $k$. But (according to Lemma \ref{generation})
      if $n$ is a generator of $G_{p^2}$,
      then $n$ is also a generator of $G_{p^l}$ for every natural number $l$ (and $p\neq 2$)
      and especially
      for $l=k$. Consequently there is an $N$ such that $n^N=x\mod p^k$. From the
      definition of the $p-$adic metric we have that $\left|x-y\right|_p<p^{-k}$ if and only if
      $x\equiv y\mod p^k$. Hence  $\left|x-n^N\right|_p<p^{-k}$. By reversing the arguments above,
      we see that the converse is also true. 
    $\square$ 

Let us consider $p\not=2$ and for $x\in B_{p^{-1}}(1)$ the $p$-adic 
exponential  function $t \mapsto x^t,$ see, for example~\cite{Schikhof}.
This function is well defined and continuous as a map from $\mathbb{Z}_p$ to
$\mathbb{Z}_p$.
In particular, for each $a \in \mathbb{Z}_p$, we have 
\begin{equation}
\label{l1}
x^a=\lim_{k\rightarrow a} x^k, \; \; k\in\mathbb{N}.
\end{equation}
We shall also use properties of the $p$-adic logarithmic function, see, for
example \cite{Schikhof}.
Let $ z \in B_{p^{-1}}(1).$  Then $\log z$ is well defined. 
For $z=1+\lambda$ with $\left|\lambda \right|_p \leq 1/p$, we have:
\begin{equation}
\label{l2}
 \log{z}=\sum_{k=1}^{\infty} \frac{(-1)^{k+1}\Delta^k}{k}=\lambda
(1+\lambda \Delta_\lambda),
      \quad \left|\Delta_\lambda \right|_p\leq 1.
\end{equation}      

\vspace{1.5ex}
\noindent  
By using (\ref{l2}) we obtain that  $\log : B_{p^{-1}}(1) \to B_{p^{-1}}(0)$ is
an isometry:
\begin{equation}
\label{l3}
\vert \log x_1 - \log x_2 \vert_p= \vert x_1 - x_2\vert_p,\; \;  x_1, x_2
\in B_{1/p}(1) \; .
\end{equation}

\begin{lemma} Let $x \in B_{p^{-1}}(1), x \not = 1, a \in \mathbb{Z}_p$ and 
let $\{ m_k \}$ be a sequence of natural numbers. If $x^{m_k} \to x^a,
k \to \infty,$ then $m_k \to a$ as $k \to \infty,$ in $\mathbb{Z}_p$.
\end{lemma}

\vspace{1.5ex}
\noindent
This is a consequence of the isometric property of $\log$.

\vspace{1.5ex}

  \begin{theorem}\label{theorem-minimal}
    Let $p\not=2$ and $l\geq 1$. Then, the monomial dynamical system
    $x\mapsto x^n$ is minimal on the circle $S_{p^{-l}}(1)$ if and only if 
$n$ is a
    generator of $G_{p^2}$.
  \end{theorem}
{\noindent \bf Proof. }%
Let $x \in S_{p^{-l}}(1).$
Consider the equation $x^a = y$.  We have that $a=\frac{\log x}{\log y}$. 
As $\log: B_{p^{-1}}(1) \to B_{p^{-1}}(0)$ is an isometry, we
have $\log(S_{p^{-l}}(1)) = S_{p^{-l}}(0).$ Thus 
$a= \frac{\log x}{\log y} \in S_1(0)$ and moreover, each
$a \in S_1(0)$ can be represented as $\frac{\log x}{\log y}$
for some $y \in S_{p^{-l}}(1)$.

Let $y$ be an arbitrary element of $S_{p^{-l}}(1)$ and 
let $x^a=y$ for some $a\in S_1(0).$ By Theorem 3.1, if $n$ is a generator of
$G_{p^2},$ 
then each $a \in S_1(0)$ is a limit point of the sequence $\{ n^N
\}_{N=1}^\infty.$
Thus $a= \lim_{k\to\infty} n^{N_k}$ for some subsequence $\{N_k \}.$ By
using the continutity
of the exponential function we obtain (\ref{limit point}).

Suppose now that, for some $n,$ $x^{n^{N_k}} \to x^a.$ By Lemma 3.2 we
obtain that
$n^{N_k} \to a$ as $k \to \infty.$ If we have (\ref{limit point}) for all $y
\in S_{p^{-l}}(1),$
then each $a\in S_1(0)$ can be approximated by elements $n^N.$ In
particular,
all elements in the set $\{1,2,...,p-1,p+1,.., p^2-1 \}$ can be approximated modulo $p^2$.
Thus $n$ is a is a generator of $G_{p^2}$. 
$\square$.

\vspace{1.5ex}
\noindent {\bf Example }In the case that $p=3$ we have that $\psi_n$ is 
minimal if $n=2$,
$2$ is a generator of $U_{3^2}=\{1,2,4,5,7,8\}$. But for $n=4$ it is not; 
$\langle 4\rangle\mod 3^2=\{1,4,7\}$. We can also see this by noting that
$S_{1/3}(1)=B_{1/9}(4)\cup B_{1/9}(7)$ and that $B_{1/9}(4)$ is invariant 
under $\psi_4$.

\begin{corollary}
If $a$ is a fixed point of the monomial dynamical system $x\mapsto x^n$, then
this is minimal on $S_{p^{-l}}(a)$ if and only if $n$ is a generator of 
$G_{p^2}$.
\end{corollary}

\vspace{1.5ex}
{\noindent \bf Proof. }%
The assertion follows immediately from Theorem 3.2 by
topological conjugation via $x\mapsto x/a$, 
$S_{p^{-l}}(a)\mapsto S_{p^{-l}}(1)$. 
$\square$

\vspace{1.5ex}
\noindent {\bf Example }
Let $p=17$ and $n=3$. On $\mathbb{Q}_{17}$ there is a primitive $3$rd root of 
unity, see for example \cite{Li1}. 
Moreover, $3$ is also a generator of $G_{17^2}$. Therefore there exist
$n$th roots of unity different from $1$ around which the dynamics is
minimal.

\section{Unique ergodicity}

In the following we will show that the minimality of the monomial dynamical
system $\psi_n:$ $x\mapsto x^n$ on the sphere $S_{p^{-l}}(1)$ is equivalent 
to its \emph{unique ergodicity}. The latter property means that there exists 
a unique probability measure on $S_{p^{-l}}(1)$ and its Borel $\sigma-$algebra
which is invariant under $\psi_n$. We will see that this 
measure is in fact the normalized restriction of the Haar measure on
$\mathbb{Z}_p$. Moreover, we will also see that the ergodicity with
respect to Haar measure of $\psi_n$ is also equivalent to its unique
ergodicity. We should point out that -- though many results are
analogous to the case of the (irrational) rotation on the circle, our
situation is quite different, in particular as we do not deal with
dynamics on topological subgroups.
\begin{lemma}\label{lemma-unique}
Assume that $\psi_n$ is minimal. Then, the Haar measure $m$ is the unique
$\psi_n-$invariant  measure on $S_{p^{-l}}(1)$.
\end{lemma}

{\noindent \bf Proof. }%
First note that minimality of $\psi_n$ implies that $(n,p)=1$ and hence that $\psi_n$
is an isometry on $S_{p^{-l}}(1)$, see equation (\ref{fundamental-inequality}).
Then, as a consequence of Theorem 27.5 in \cite{Schikhof}, it follows that 
$\psi_n(B_r(a))=B_r(\psi_n(a))$ for each ball $B_r(a)\subset S_{p^{-l}}(1)$.
Consequently, for every open set $U\neq\emptyset$ we must have that 
$S_{p^{-l}}(1)=\cup_{N=0}^{\infty}\psi_n^N(U)$. 
Then, it follows for a $\psi_n-$invariant measure $\mu$ that $\mu (U)>0$. 

Moreover, we can split $S_{p^{-l}}(1)$ into disjoint balls of radii
$p^{-(l+k)}$, $k\geq 1$, on which $\psi_n$ acts as a permutation. 
In fact, for each $k\geq 1$, $S_{p^{-l}}(1)$ is the union,
\begin{equation}\label{splitting}
S_{p^{-l}}(1)=\cup B_{p^{-(l+k)}}(1+b_lp^l+...+b_{l+k-1}p^{l+k-1}),
\end{equation}
where $b_i\in\{0,1,...,p-1\}$ and $b_l\neq 0$.

We now show that $\psi_n$ is a (transitive) permutation on the partition
(\ref{splitting}). Recall that  every element of a $p-$adic ball is
the center of that ball (see section \ref{p-adic numbers}), and as
pointed out above $\psi_n(B_r(a))=B_r(\psi_n(a))$. Consequently we
have for all positve integers $k$, $\psi_n^k(a)\in B_r(a)\Rightarrow \psi_n^k(B_r(a))=B_r(\psi_n^k(a))=B_r(a)$ so that
$\psi_n^{Nk}(a)\in B_r(a)$ for every natural number
$N$. Hence, for a minimal $\psi_n$ a point of a ball $B$ of the
partition (\ref{splitting}) must move to another ball in the partition.
Furthermore, the minimality of $\psi_n$
shows indeed that $\psi_n$ acts as
a permutation on balls. By invariance of $\mu$, all balls must have the 
same positive measure. As this holds for any $k$, $\mu$ must be the
restriction of Haar measure $m$. 
$\square$

\vspace{1.5ex}
\noindent
The arguments of the proof of Lemma \ref{lemma-unique} also show that Haar 
measure is always $\psi_n-$invariant. Thus if $\psi_n$ is uniquely ergodic, 
the unique invariant measure must be the Haar measure $m$. Under these 
circumstances it is known (\cite[Theorem 6.17]{Walters}) that $\psi_n$ must 
be minimal.

\begin{theorem}\label{theorem-minimal-unique}
The monomial dynamical system $\psi_n:$ $x\mapsto x^n$ on $S_{p^{-l}}(1)$ is 
minimal if and only if it is uniquely ergodic in which case the unique 
invariant measure is the Haar measure.
\end{theorem}

\vspace{1.5ex}
\noindent
Let us mention that unique ergodicity yields in particular the ergodicity of 
the unique invariant measure, i.e. the Haar measure $m$, which means that
\begin{equation}
\frac {1}{N}\sum_{i=0}^{N-1}f(x^{n^i})\rightarrow\int\! f\, dm \text{ for all }
x\in S_{p^{-l}}(1),
\end{equation}
and all continous functions $f:\quad S_{p^{-l}}(1)\rightarrow \mathbb{R}$.

On the other hand the arguments of the proof of Lemma \ref{lemma-unique}, i.e.
 the fact that $\psi_n$ acts as a transitive permutation on each partition of 
$S_{p^{-l}}(1)$ into disjoint balls if and only if 
$\langle n\rangle =G_{p^2}$,  proves that if $n$ is not a generator of 
$G_{p^2}$ then the system is not ergodic with respect to Haar measure. Consequently, 
if $\psi_n$ is ergodic then $\langle n\rangle =G_{p^2}$ so that the system is 
minimal by Theorem \ref{theorem-minimal},
and hence even uniquely ergodic by Theorem \ref{theorem-minimal-unique}. Since
 unique ergodicity implies ergodicity  one has the following.

\begin{theorem}
The monomial dynamical system $\psi_n:$ $x\mapsto x^n$ on $S_{p^{-l}}(1)$ is 
ergodic with respect to Haar measure if and only if it is uniquely ergodic.
\end{theorem}

Even if the monomial dynamical system $\psi_n:$ $x\mapsto x^n$ on
$S_{p^{-l}}(1)$ is ergodic, it cannot be mixing, especially not
weak-mixing. This can be seen from the fact that an abstract dynamical
system is weak-mixing if and only if the product of such two systems
is ergodic. If we choose a function $f$ on $S_{p^{-l}}(1)$ and define
a function $F$ on $S_{p^{-l}}(1)\times S_{p^{-l}}(1)$ by
$F(x,y):=f(\log x/\log y)$ (which is well defined as $\log$ does not
vanish on $S_{p^{-l}}(1)$), we obtain a non-constant function
satisfying $F(\psi_n(x),\psi_n(y))=F(x,y)$. This shows (see
\cite[Theorem 1.6]{Walters}) that $\psi_n\times\psi_n$ is not ergodic,
and hence $\psi_n$ is not weak-mixing with respect to any invariant
measure, in particular not to the Haar measure.

\vspace{1.5ex}
\noindent
Let us consider the ergodicity of a perturbed system
\begin{equation}\label{perturbation}
\psi_q=x^n+q(x),
\end{equation}
for some polynomial $q$ such that $q(x)\equiv 0\mod p^{l+2}$, 
($\left | q(x)\right |_p\leq p^{-(l+2)}$). Under this  condition, the sphere $S_{p^{-l}}(1)$ is invariant. For such a system to be ergodic
 it is necessary that $n$ is a generator of $G_{p^2}$. This follows from the
 fact that for each $x=1+a_lp^l+a_{l+1}p^{l+1}+...$  on $S_{p^{-l}}(1)$ 
 the condition on $q$ gives 
\begin{equation}
\psi_q^N(x)\equiv 1+n^N(a_l+a_{l+1}p)p^l \mod p^{l+2}.
\end{equation}
Now, $\psi_q$ acts as a transitive permutation on the $p(p-1)$ balls of radius 
$p^{-(l+2)}$ if and only if $\langle n\rangle =G_{p^2}$. Consequently, a 
perturbation (\ref{perturbation}) cannot make a nonergodic system ergodic.
However, it is not clear at this point if this is true also if $|q(x)|=p^{-(l+1)}$.  
  
\section{Acknowledgement}
We would like to thank P.A. Svensson at V\"axj\"o Univerity for consultation in algebra.

\appendix

\section{Appendix}\label{lemma}

\begin{lemma}\footnote{This result was obtained in \cite{Li1}.}
Let $x,y\in S_1(0)$ and suppose $\left| x-y\right| _p<1.$
Then for all natural numbers $n,$%
\begin{equation}
\left| x^n-y^n\right| _p\leq \left| n\right| _p\left| x-y\right| _p,
\end{equation}
with equality for $p>2$. Moreover equality also holds for $p=2$ if $n$
is odd.
\end{lemma}

{\noindent \bf Proof. }First note that 
\begin{eqnarray*}
\left| x^n-y^n\right| _p &=&\left| (y+x-y)^n-y^n\right| _p=\left|
\sum_{k=0}^n\binom nky^{n-k}(x-y)^k-y^n\right| _p \\
&=&\left| \sum_{k=1}^n\binom nky^{n-k}(x-y)^k\right| _p\leq \max_k\left|
\binom nky^{n-k}(x-y)^k\right| _p,
\end{eqnarray*}
where the inequality is a consequence of the strong triangle inequality (\ref{str}). We
will see that the maximum occurs for $k=1$ and that the corresponding term is
greater than every other term (except for the case where $p=2$ and $n$ is even).
By the multiplicative property
of the valuation, the valuation of the $k$th term in the sum can be
estimated according to 
\begin{align*}
\left| \binom nky^{n-k}(x-y)^k\right| _p &=\left| \binom nk\right| _p\left|
y\right| _p^{n-k}\left| x-y\right| _p^k \\
&=\left| n(n-1)\cdot ...\cdot (n-k+1)\right| _p\left| x-y\right| _p\frac{%
\left| x-y\right| _p^{k-1}}{\left| k!\right| _p} \\
&\leq \left| n\right| _p\left| x-y\right| _p\frac{\left| x-y\right| _p^{k-1}%
}{\left| k!\right| _p}\\&=\left| \binom n1y^{n-1}(x-y)^1\right| _p\frac{\left|
x-y\right| _p^{k-1}}{\left| k!\right| _p},
\end{align*}
where the inequality holds because $\left| a\right| _p\leq 1$ for every
integer $a\in \mathbb{Z}$ (together with the multiplicative property of the
valuation). Note that strong inequality holds in the case that $p=2$ and $n$ is odd;
if $n$ is odd we have $\left| n\right| _2=1$ ($n$ is not divisible by $2$). Moreover
we can assume that $n>1$ (the lemma is automatically true for $n=1$) so that 
$\left| n-1\right| _2<1$.

To estimate the latter expression, $\left| x-y\right| _p^{k-1}/\left| k!\right| _p$,
we use the fact that
$\left| k!\right| _p\geq 1/p^{k-1}$ with strong inequality for $p>2$,
see for example the paper \cite[Lemma 25.5]{Schikhof}.
Moreover $\left| x-y\right| _p\leq 1/p$ by the condition of the lemma.
Consequently, 
\[
\left| x-y\right| _p^{k-1}\leq 1/p^{k-1}. 
\]
Hence we have that $\left| x-y\right| _p^{k-1}/\left| k!\right| _p\leq 1$
with strong inequality when $p>2.$ Therefore, for every integer $k$ such
that $1\leq k\leq n,$ we have 
\[
\left| \binom nky^{n-k}(x-y)^k\right| _p\leq \left| \binom
n1y^{n-1}(x-y)^1\right| _p, 
\]
with strong inequality for all $p$ and $n$ except the case  when $p=2$ and $n$ is even.
Hence, as a consequence of the strong triangle inequality (\ref{str}) we have that 
\[
\left| x^n-y^n\right| _p\leq \left| \binom n1y^{n-1}(x-y)^1\right| _p=\left|
n\right| _p\left| x-y\right| _p, 
\]
with equality for $p>2$ and, if $n$ is odd, also for $p=2$.
$\Box $

\end{document}